\documentclass[a4paper]{mathscan}

\usepackage{epsfig,color}

\usepackage{mathrsfs,enumerate}
\newtheorem{theorem}{Theorem}[section] 
\newtheorem{prop}[theorem]{Proposition} 
\newtheorem{cor}[theorem]{Corollary}    
\newtheorem{lemma}[theorem]{Lemma}        
             \newtheorem*{thma}{Theorem}

\renewcommand{\mathbb}{\mathsf}

\def\codim{\operatorname{codim}}

\def\max{\operatorname{max}}

\def\c1{\operatorname{c_1}}
\def\c2{\operatorname{c_2}}

\def\Cox{\operatorname{Cox}}

\def\ZZ{{\mathbf Z}}

\def\PP{{\mathbf P}}

\def\L{{\mathscr L}}

\def\O{{\mathscr O}}

\def\F{{\mathscr F}}

\def\+{\oplus}                   
\def\*{\otimes}                  

\def\Pic{\operatorname{Pic}}

\renewcommand{\qed}{\hfill \ensuremath{\Box}}

\title[The Cox Ring of $\PP^2$ Blown Up in Points on a Line]{On the Cox Ring of $\PP^2$ Blown Up in Points on a Line}
\author{John Christian Ottem}
\date{\today}

\address{Matematisk institutt, Universitetet i Oslo, PO Box 1053, Blindern, NO-0316 Oslo, Norway}
\email{johnco@math.uio.no}

\begin{document}

\begin{abstract} 
We show that the blow-up of $\PP^2$ in $n$ points on a line has finitely generated Cox ring. We give explicit generators for the ring and calculate its defining ideal of relations. 
\end{abstract}

\maketitle

\section{Introduction}In 2000, Hu and Keel \cite{HK00} introduced the \emph{Cox ring} of an algebraic variety, aiming to generalize the Cox construction for toric varieties. If $X$ is a normal variety with freely finitely  generated Picard group $\Pic(X)$, this ring is essentially defined by
\begin{equation*}
\Cox(X)=\bigoplus_{D\in \Pic(X)}H^0(X,\O_X(D)),
\end{equation*}
where the ring product is given by multiplication of sections as rational functions. Varieties whose Cox ring is finitely generated are called \emph{Mori dream spaces}, and have interesting properties from the viewpoint of birational geometry. 

Even though the definition of the Cox ring is quite explicit, calculating its presentation for a concrete variety can be a hard problem. For example, the Cox rings of Del Pezzo surfaces have been the subject of much recent literature in algebraic geometry (e.g., \cite{BP04}, \cite{LV07} and \cite{STV06}) and show that the behaviour of the Cox ring under blow-up is highly non-trivial. For general blow-ups of $\PP^2$, $\Cox(X)$ may even fail to be finitely generated, since the surfaces may have infinitely many curves of negative self-intersection.

In this paper, we consider blow-ups of $\PP^2$ in points $p_1,\ldots, p_n$ lying on a fixed line $Y \subset \PP^2$. The blow-up $\pi: X\to \PP^2$ in these points is a smooth projective surface with Picard group generated by the classes of the exceptional divisors $E_1,\ldots, E_n$ and  $L$, the pullback of a general line in $\PP^2$. Our first result is that $\Cox(X)$ is finitely generated for any number of points on a line. This result was first shown in \cite{EKW04} and correlates with recent results of Hausen and S\"uss in \cite{HS09}, since the surface $X$ has a complexity one torus action.

Furthermore, we also find explicit generators for $\Cox(X)$, i.e., generating sections $x_1,\ldots,x_r$ from respective vector spaces $H^0(X,\O_X(D_1)),\ldots, H^0(X,\O_X(D_r))$, so that $\Cox(X)$ may be regarded as a quotient
$$
\Cox(X)\simeq k[x_1,\ldots, x_r]/I.
$$Here we consider a $\Pic(X)-$grading on $k[x_1,\ldots, x_r]$ and $I$ given by $\deg(x_i)=D_i$. In Section 4 we find explicit generators and a Gr\"obner basis for the ideal $I$. Our main result is the following theorem:

\begin{thma}Let $X$ be the blow-up of $\PP^2$ in $n\ge 3$ distinct points lying on a line. Then $\Cox(X)$ is a complete intersection ring and its defining ideal is generated by quadric trinomials.
\end{thma}

\noindent In particular, this means that the Cox ring is Gorenstein and a Koszul algebra.

\bigskip 

\noindent \emph{Notation:} We make the following standard shorthand notation for sheaf cohomology:
$$
H^i(D):=H^i(X,\O_{X}(D)), \qquad  h^i(D):=\dim_k H^i(X,\O_{X}(D)), \qquad i=0,1,2.
$$ 

\section{Nef divisors and vanishing on $X$. }
The surface $X$ has good vanishing properties for nef divisors. For example, $H^2(D)=H^0(K-D)=0$ is immediate by Serre duality, since $K$ cannot be effective on $X$. It turns out that also $H^1(D)=0$ for $D$ nef, so all cohomology can be calculated from the Riemann-Roch theorem. To prove this, we first need some preparatory lemmas.

\begin{lemma}\label{eff}
The monoid of effective divisor classes of $X$ is finitely generated by the classes $L-E_1-\ldots-E_n,E_1,E_2,\ldots,E_n$.
\end{lemma}\begin{proof}
It is clear that the classes above are all effective, so their semigroup span is in the effective monoid. Conversely, note that these divisor classes actually form a $\ZZ$-basis for $\Pic X$. So given an \emph{irreducible} effective divisor $D$, we let
$$
m(L-E_1-\ldots-E_n)+a_{1}E_1+\ldots+a_nE_n
$$represent the corresponding divisor class. If $D$ is not one of the generators above we have $D.E_i=m-a_i\ge 0$ and $D.(L-E_1-\ldots-E_n)=m-\sum_{i=1}^n(m-a_i)\ge 0$. Together these inequalities imply that $m,a_i\ge 0$, and we are done.\end{proof}
\begin{lemma}\label{nefcone}
The nef monoid is generated by the divisor classes $L,L-E_1,L-E_2,\ldots, L-E_n$.
\end{lemma}

\begin{proof} Let $D=dL-\sum a_iE_i$ be a nef divisor class. Intersecting $D$ with the classes in Lemma \ref{eff} gives the following set of inequalities:
$$
d\ge a_1+a_2+\ldots+a_n,\qquad a_i\ge 0, \quad \forall i=1,\ldots,n
$$Now it is easy to see that we can decompose each $D$ as a sum of the $L-E_i$'s by using $a_i$ of $L-E_i$ and finally add $d-a_1-a_2-\ldots-a_n\ge 0$ times $L$. 
\end{proof}\noindent Note that since the classes above are effective, so is every nef divisor on $X$.

\begin{lemma}\label{h1}
Let $D=dL-a_1E_1-\ldots-a_nE_n$ be a divisor class on $X$, with $d+1\ge \sum_{i=1}^n a_i$ and $a_i\ge 0$. Then $h^1(D)=0$.
\end{lemma}\begin{proof}If $D=dL$, we have $h^1(D)=h^1(\PP^2,\O_{\PP^2}(d))=0$ for $d\ge -1$. If say $a_1>0$, consider the divisor class $D'=D-(L-E_1)$. $D'$ satisfies the conditions of the lemma, so by induction on $d$ we have $h^1(D')=0$. Let $C$ be a smooth rational curve with class $L-E_1$, then $h^1(C,D|_C)=h^1(\PP^1,\O_{\PP^1}(D.C))=0$, since $D.C\ge -1$. Now taking cohomology of the exact sequence
$$
0\to \O_X(D')\to \O_X(D)\to \O_C(D|_C)\to 0
$$gives $h^1(D)=0$. 
\end{proof}

\noindent This gives us the multigraded Hilbert function of $\Cox(X)$ in nef degrees:

\begin{cor}\label{hilbert}
For a nef divisor class $D$, we have $h^1(D)=0$ and $$\dim_k \Cox(X)_D = \chi(\O_X(D))$$.
\end{cor}

\section{Generators for $\Cox(X)$.}
We are now in position to find explicit generators for $\Cox(X)$ as a $k-$algebra. 
When $n=1$, the surface $Bl_{p}\PP^2$ is a toric variety, and by \cite{Cox95}, its Cox ring coincides with the 
usual homogenous coordinate ring
$$
\Cox( Bl_{p}\PP^2)=k[x,s_1,s_2,e] 
$$where $\deg x=L$, $\deg s_i=L-E$ and $\deg e=E$. Therefore, we will in the following suppose that $n\ge 2$.

We first choose generators $e_1,\ldots,e_n$ for the 1-dimensional vector spaces $H^0(E_i)$ for $i=1,\ldots,n$ and a generator $l$ of $H^0(L-E_1-\ldots-E_n)$. For the classes $L-E_i$, for which $H^0(L-E_i)$ is $2$-dimensional, we need in addition to the section $le_1\cdots e_{i-1}e_{i+1}\cdots e_n$, a new section $s_i$ to form a basis. To specify these explicitly, we fix a point $q\in \PP^2$ and for each $i$ take a section corresponding to the strict transform of the line going through $q$ and $p_i$. The projections of these sections to $\mathbb{P}^2$ are shown in Figure \ref{3lines}:

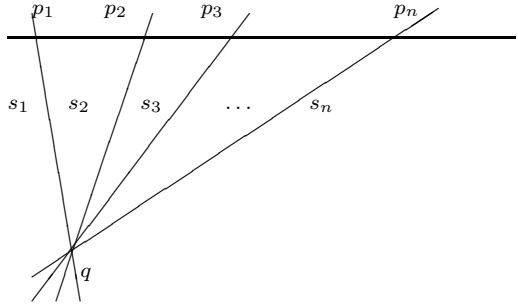
\begin{figure}[htbp]
\begin{center}
\setlength{\unitlength}{3947sp}%
\begingroup\makeatletter\ifx\SetFigFont\undefined%
\gdef\SetFigFont#1#2#3#4#5{%
  \reset@font\fontsize{#1}{#2pt}%
  \fontfamily{#3}\fontseries{#4}\fontshape{#5}%
  \selectfont}%
\fi\endgroup%
\begin{picture}(3174,1908)(889,-2323)
\thinlines
{\color[rgb]{0,0,0}\put(901,-661){\line( 1, 0){3150}}
}%
{\color[rgb]{0,0,0}\put(1051,-511){\line( 1,-6){300}}
}%
{\color[rgb]{0,0,0}\put(2401,-511){\line(-3,-4){1350}}
}%
{\color[rgb]{0,0,0}\put(1201,-2311){\line( 1, 3){600}}
}%
{\color[rgb]{0,0,0}\put(1051,-2161){\line( 3, 2){2526.923}}
}%
\put(1351,-2161){\makebox(0,0)[lb]{\smash{{\SetFigFont{8}{9.6}{\familydefault}{\mddefault}{\updefault}{\color[rgb]{0,0,0}$q$}%
}}}}
\put(1051,-511){\makebox(0,0)[lb]{\smash{{\SetFigFont{8}{9.6}{\familydefault}{\mddefault}{\updefault}{\color[rgb]{0,0,0}$p_1$}%
}}}}
\put(1501,-511){\makebox(0,0)[lb]{\smash{{\SetFigFont{8}{9.6}{\familydefault}{\mddefault}{\updefault}{\color[rgb]{0,0,0}$p_2$}%
}}}}
\put(2101,-511){\makebox(0,0)[lb]{\smash{{\SetFigFont{8}{9.6}{\familydefault}{\mddefault}{\updefault}{\color[rgb]{0,0,0}$p_3$}%
}}}}
\put(3301,-511){\makebox(0,0)[lb]{\smash{{\SetFigFont{8}{9.6}{\familydefault}{\mddefault}{\updefault}{\color[rgb]{0,0,0}$p_n$}%
}}}}
\put(901,-1111){\makebox(0,0)[lb]{\smash{{\SetFigFont{8}{9.6}{\familydefault}{\mddefault}{\updefault}{\color[rgb]{0,0,0}$s_1$}%
}}}}
\put(1276,-1111){\makebox(0,0)[lb]{\smash{{\SetFigFont{8}{9.6}{\familydefault}{\mddefault}{\updefault}{\color[rgb]{0,0,0}$s_2$}%
}}}}
\put(1726,-1111){\makebox(0,0)[lb]{\smash{{\SetFigFont{8}{9.6}{\familydefault}{\mddefault}{\updefault}{\color[rgb]{0,0,0}$s_3$}%
}}}}
\put(2776,-1111){\makebox(0,0)[lb]{\smash{{\SetFigFont{8}{9.6}{\familydefault}{\mddefault}{\updefault}{\color[rgb]{0,0,0}$s_n$}%
}}}}
\put(2251,-1111){\makebox(0,0)[lb]{\smash{{\SetFigFont{8}{9.6}{\familydefault}{\mddefault}{\updefault}{\color[rgb]{0,0,0}$\ldots$}%
}}}}
\end{picture}%

\caption{The choice of the sections $s_1,s_2,\ldots, s_n$.}
\label{3lines}
\end{center}
\end{figure}

\noindent We will show that $\Cox(X)$ is generated by the sections $l,e_i,s_i$ for $i=1,\ldots,n$. The following lemma is a variant of Castelnuovo's base point free pencil trick (\cite[Ex. 17.18]{Eis95}) and will be our main tool for proving this.

\begin{lemma}\label{bfpt} Let $X$ be an algebraic variety over a field $k$, let $\F$ be a locally free sheaf of $\O_X$-modules on $X$, let  $\L$ be an invertible sheaf on $X$ and $V$ a two-dimensional base-point free subspace of $H^0(X,\L)$. If $H^1(\L^{-1}\otimes \F)=0$, then the multiplication map
\begin{equation*}
V\otimes H^0(X,\F)\to H^0(X,\L\otimes \F)
\end{equation*}is surjective.
\end{lemma}

\begin{prop}\label{linegens}
Let $X$ be the blow-up of $\PP^2$ in $n\ge 2$ distinct points on a line. Then there is a multigraded surjection
\begin{equation}\label{surjection}
k[l,e_1,\ldots, e_n,s_1,\ldots,s_n]\to \Cox(X).
\end{equation}where $\deg(l)=L-E_1-\ldots-E_n$, $\deg e_i=E_i$ and $\deg s_i=L-E_i$.
\end{prop}

\begin{proof}Let $D$ be an effective divisor on $X$. We need to show that $H^0(D)$ has a basis of sections which are polynomials in $l, e_i, s_i$.

We first show that we may take $D$ to be nef. Indeed, suppose $E$ is a curve such that $D\cdot E<0$. Without loss of generality, we may suppose that $E$ is one of the divisor classes generating the effective monoid. Let $x_E\in \{l,e_1,\ldots,e_n\}$ be the corresponding section in $H^0(E)$. Then since $E$ is a base component of the linear system $|D|$, multiplication by $x_E$ induces an isomorphism $H^0(D-E)\to H^0(X,D)$. By induction on the number of fixed components, $H^0(D-E)$ is generated by monomials in the $l,e_i,s_i$ and hence the same applies to $H^0(D)$.

Now suppose that $D$ is a nef divisor class and write $D$ (uniquely) in terms of the nef classes of Lemma \ref{nefcone}: 
$$
D=aL+a_1(L-E_1)+a_2(L-E_2)+\ldots+a_n(L-E_n)
$$where $a, a_i\ge 0$. 
If say, $a_1\ge 2$, then $H^1(D-2(L-E_1))=0$, since $D-2(L-E_1)$ is nef, and so applying Lemma \ref{bfpt} with $V=H^0(L-E_1)$ and $\F=\O_X(D-(L-E_1))$, we get a surjection
$$
H^0(D-(L-E_1))\otimes H^0(L-E_1)\to H^0(D).
$$By induction on the number $D.L\ge0$, $H^0(D-(L-E_1))$ is generated by monomials in $l,e_i,s_i$, and therefore so is $H^0(D)$. 

If $a_i\le 1$ for all $i$, and say, $a_1=1$, then $D-2(L-E_1)=N+E_1$ for some divisor $N$ satisfying the assumptions of Lemma \ref{h1} and $N.E_1=0$. In particular, $h^1(N)=0$. Now also $h^1(N+E_1)=0$, by the exact sequence
$$
0\to \O_{X}(N)\to \O_X(N+E_1) \to \O_{E_1}(-1)\to 0
$$ and we proceed as above. 

If $a_i=0$ for all $i$, then $D=aL$  for some $a\ge 1$ and $H^0(X,\O_X(D))\simeq H^0(\PP^2,\O_{\PP^2}(a))$. This implies that $$H^0(X,(a-1)L)\otimes H^0(X,L)\to  H^0(X,aL).$$is surjective. By induction on $a$, $H^0((a-1)L)$ is generated by monomials in $l,e_i,s_i$, and therefore so is $H^0(D)$.

It remains to show that $H^0(L)$ has a basis of monomials in $l,e_i,s_i$. But $H^0(L)=\pi^*H^0(\PP^2,\O_{\PP^2}(1))$, so it suffices to find three monomials of degree $L$ that project to linearly independent sections in $\O_{\PP^2}(1)$. By construction, this works for the three sections
$$
\sigma_1=s_1e_1,\quad\sigma_2 =s_2e_2, \quad\sigma_3=  le_1e_2e_3\cdots e_n.
$$
\end{proof}

\section{Relations}
\noindent We now turn to the defining ideal $I$ of relations of $\Cox(X)$, i.e., the kernel of the map (\ref{surjection}). 
Consider again the divisor class $L$: We have $h^0(L)=3$, while there are $n+1$ monomials of degree $L$ in $k[l,e_i,s_i]$  :
$$
s_1e_1,\quad s_2e_2, \quad\cdots \quad s_ne_n, \quad le_1e_2e_3\cdots e_n
$$This means that there are $n-2$ linear dependence relations between them. To see what they look like, consider again the projection of these sections in Figure \ref{3lines}. Of course any three of these lines through $q$ satisfy a linear dependence relation, and these pull back via $\pi$ to relations in $\Cox(X)$ of the following form:\begin{eqnarray}\label{eqns}
g_1&=&\underline{s_1e_1}+a_1s_{n-1}e_{n-1}+b_1s_ne_n=0\nonumber \\
g_2&=&\underline{s_2e_2}+a_2s_{n-1}e_{n-1}+b_2s_ne_n=0 \nonumber\\
\vdots & & \qquad \qquad \qquad \qquad \vdots\\
g_{n-2}&=&\underline{s_{n-2}e_{n-2}}+a_{n-2}s_{n-1}e_{n-1}+b_{n-2}s_ne_n=0 \nonumber
\end{eqnarray}where each of the coefficients $a_i,b_i$ are non-zero. We denote the ideal generated by these relations by $J$. The leftmost terms above are underlined since as the next lemma shows, they form an initial ideal for $J$. 

\begin{lemma}
The set $\{g_1,\ldots,g_{n-2}\}$ is a Gr\"obner basis for $J$ with respect to the graded lexicographical order, and $(s_1e_1,\ldots,s_{n-2}e_{n-2})$  is an initial ideal of $J$.
\end{lemma}\begin{proof} It is well-known (e.g., see \cite{AL94}) that a collection of polynomials with relatively prime leading terms is a Gr\"obner basis for the ideal they generate.  \end{proof}

\noindent We will show that that the expressions (\ref{eqns}) in fact generate all the relations, i.e., that $I=J$. For this, we will make use of the $\Pic(X)-$grading on $R=k[l,e_i,s_i]$ and $I$. The next lemma shows that it is sufficient to consider generators for the ideal of degrees corresponding to nef divisor classes.

\begin{lemma}\label{idealnef}
The ideal $I$ is generated by elements of degree $D$, where $D$ is a nef divisor class.
\end{lemma}\begin{proof}Suppose $D$ is an effective divisor class and that there is a negative curve $E$ such that $D.E< 0$. Then this implies that $E$ is a fixed component of $|D|$ and as above every monomial in $k[l,s_i,e_i]_D$ is divisible by $x_E$, the variable corresponding to $E$. This shows that any element of $I_D$ can be written as a product of $x_E$ and a relation in $I_{D-E}$. Now the claim follows by induction on the number of fixed components of $D$.  \end{proof}

We will now prove our main theorem.

\begin{theorem}
Let $X$ be the blow-up of $n\ge 2$ points on a line. Then $\Cox(X)$ is a complete intersection with $n-2$ quadratic defining relations given in (\ref{eqns}).
\end{theorem}
\begin{proof}By Lemma \ref{idealnef}, it is sufficient to show that $I_{D}=J_{D}$ for all nef classes $D=dL-a_1E_1-a_2E_2-\ldots-a_nE_n$, (here $d\ge a_1+\ldots+a_n$). Note that since $J\subseteq I$, we have in any case a surjective homomorphism
$$
R/J\to \Cox(X).
$$To show that this is an isomorphism in degree $D$, we calculate the multigraded Hilbert function of both sides. From Corollary \ref{hilbert} and Riemann-Roch, we have 
\begin{equation}\label{RRabc}
\dim_k \Cox(X)_{D}=h^0(D)={d+2\choose 2}-{a_1+1\choose 2}-\ldots-{a_n+1\choose 2}.
\end{equation}
To calculate $\dim_k ({R/J})_{D}$, we use the Gr\"obner basis for $J$. Since the Hilbert function is preserved when going to initial ideals, we have 
$$
\dim_k (R/J)_{D}=\dim_k {R/(s_1e_1,\ldots,s_{n-2}e_{n-2})}_{D}
$$We use a counting argument to evaluate this number. Note that any monomial in ${R/(s_1e_1,\ldots,s_{n-2}e_{n-2})}_{D}$ corresponds to a way of writing $D$ as a non-negative sum of divisor classes from $$L-E_1-\ldots-E_n, \qquad L-E_1,\qquad \ldots \qquad L-E_n,\qquad E_1, \qquad\ldots, \qquad E_n$$such that not both $L-E_i$ and $E_i$ occur in the sum for $i=1,\ldots,n-2$.  Abusing notation, we let the numbers $s_i ,e_i,l$ represents respectively the non-negative coefficients of $L-E_i,E_i,L-E_1-\ldots-E_n$ in this sum. Working in $\Pic(X)\simeq \mathbb{Z}^{n+1}$, this translates the problem of finding $\dim_k (R/J)_D$ into following counting problem: finding the number of non-negative solutions of
\begin{eqnarray}\label{mabc'}
s_1+s_2+\ldots+s_n &=& d-l\nonumber\\
s_1-e_1 &=& a_1-l\nonumber\\
\vdots & &\vdots\nonumber\\
s_n-e_n &=& a_n-l\nonumber
\end{eqnarray}such that $s_i\cdot e_i=0$ for $i=1,\ldots,n-2.$ 

\begin{lemma} For each $l \le d$ we have $d-l\ge \sum_{k=1}^n\max(a_i-l,0)$.
\end{lemma}\begin{proof}For $l=0$, this inequality reduces to the nef condition on $D$. Now, increasing $l$ by one decreases the left hand side by one, and if there is some $a_i-l>0$, then $\max(a_i-l,0)$ is decreased by one, as well, if not, the right hand side is zero, so in any case the inequality is preserved.
\end{proof}
For each fixed $l$, we count the number of non-negative solutions $S(l)$ to the system (\ref{mabc'}). Note that $s_i$  is completely determined as $s_i=\max(a_i-l,0)$ for $1\le i \le n-2$ by the condition $s_i.e_i=0$. Hence by the first equation in (\ref{mabc'}), we are looking for non-negative solutions to 
$$
s_n+s_{n-1}=d-l-\sum_{k=1}^{n-2}\max(a_k-l,0)\ge 0
$$such that $s_n\ge \max(a_n-l,0)$ and $s_{n-1}\ge \max(a_{n-1}-l,0)$, of which there are in total
$$d-l-\sum_{k=1}^{n-2}\max(a_{n-2}-l,0)+1-\sum_{k=n-1,n}\max(a_k-l,0)$$
Hence the total number of solutions to (\ref{mabc'}) is \begin{equation*}\begin{split}
\sum_{l=0}^d S(l)&=\sum_{l=0}^d \left(d+1-l-\sum_{k=1}^n\max(a_k-l,0)\right)\\
&={d+2\choose 2}-\sum_{i=0}^{a_1} (a_1-i) - \sum_{i=0}^{a_1} (a_1-i) - \ldots - \sum_{i=0}^{a_n} (a_n-i)  \\
&={d+2\choose 2}-{a_1+1\choose 2}-{a_2+1\choose 2}-\ldots-{a_n+1\choose 2}=h^0(D).
\end{split}\end{equation*}This finishes the proof that $I=J$. Now, from \cite[Remark 1.4]{BP04} we have $\dim \Cox(X)=n+3$, furthermore by Proposition \ref{linegens} we have that $\codim \Cox(X)=(2n+1)-(n+3)=n-2$, which is exactly the number of relations in $I$.
\end{proof}

\emph{Remark. }The above result can also be proved in another way, using the following lemma, proved by Stillman in \cite{ST05}:

\begin{lemma}
Let $J \subset k[x_1, x_2, \ldots , x_n]$ be an ideal containing a polynomial
$f = gx_1 + h$, with $g$, $h$ not involving $x_1$ and $g$ a non-zero divisor modulo $J$. Then, $J$ is prime if and only if the elimination ideal $J\cap k[x_2, \ldots,x_n]$ is prime.
\end{lemma} The above lemma can be used to prove that $J=(g_1,\ldots,g_{n-2})$ is prime, using induction on $n$. For $n=3$, this is clear. Next, note that the elimination ideal $J\cap k[s_2,\ldots,s_n,e_2\ldots,e_n,l]$ is just $(g_2,\ldots,g_{n-2})$ since $\{g_1,\ldots,g_{n-2}\}$ is a Gr\"obner basis. By induction on $n$, $(g_2,\ldots,g_{n-2})$ is the defining ideal of $\Cox(X)$ for $\PP^2$ blown up in the points $p_2,\ldots,p_n$ and hence is prime.  Taking now $x_1=e_1, g=s_1$, $h=a_1s_{n-1}e_{n-1}+b_1s_ne_n$ shows that $J$ is prime. Then, since $J \subseteq I$ are two prime ideals with the same Krull dimension, they must be equal.

\bigskip

\emph{Acknowledgement.} I wish to thank Kristian Ranestad, J\o rgen Vold Rennemo and Mauricio Velasco for interesting discussions. I would also like to thank Brian Harbourne for pointing out the reference \cite{GM05} and the anonymous referee for helpful comments.

\end{document}